\documentclass[reqno,11pt]{amsart}

\usepackage[sorted-cites,abbrev]{amsrefs} 
\usepackage{amsfonts}
\usepackage{amsmath,amsthm,amscd, amssymb} 
\usepackage{amsthm}

\usepackage{color}

\title[Bergman kernel expansion]{Abstract Bergman kernel expansion and its applications}
\author{Chiung-ju Liu}
\address[Chiung-ju Liu]{Department of Mathematics, National Taiwan
University, Taipei, Taiwan 106} 
\email{cjliu4@ntu.edu.tw}
\author{Zhiqin Lu}
\address[Zhiqin Lu]{Department of
Mathematics, University of California, Irvine, CA 92697-3875}
\date{March 19, 2014}


\thanks{The first author is supported by the NSC grant NSC 982115M002007
in Taiwan. The second author is partially supported by the   NSF grant DMS-12-06748.} 
\email{zlu@uci.edu}\newtheorem{theorem}{Theorem}[section]

\newtheorem{lemma}{Lemma}[section]
\newtheorem{cor}{Corollary}[section]
\newtheorem{prop}{Proposition}[section]

\newtheorem{definition}{Definition}[section]

\theoremstyle{remark}
\newtheorem{rem}{Remark}[section]

\renewcommand{\bar}{\overline}
\newcommand{\eps}{\epsilon}
\newcommand{\pa}{\partial}

\newcommand{\rr}{{|z| \leq\frac{\log m}{\sqrt m}}}
\newfont{\fnt}{cmr10 scaled 550}
\renewcommand{\eps}{\varepsilon}

\newcommand{\ka}{K\"ahler }

\newcommand{\C}{\mathbb C}

\newcommand{\Z}{\mathbb Z}

\newcommand{\bb}{{\frac{\sqrt{-1}}{2\pi}}}

\usepackage{color}

\renewcommand{\phi}{\varphi}

\usepackage{amsmath,amsthm,amscd}

\begin{document}

\begin{abstract} We give a purely complex geometric proof of the existence of the Bergman kernel expansion. Our method actually provides a sharper estimate, and in the case that the metrics are real analytic, we prove that the remainder decays faster than any polynomial.
\end{abstract}

\maketitle \tableofcontents\pagestyle{myheadings}

\section{Introduction}

Let $M$ be an $n$-dimensional algebraic manifold in a
certain projective space $\mathbb{CP}^N$. The hyperplane line
bundle of $\mathbb{CP}^N$,  restricting to an ample line bundle $L$ of
$M$,  is  a polarization of $M$. A K\"ahler metric $g$ is
called a polarized metric, if the corresponding K\"ahler form
represents the first Chern class $c_1(L)$ of $L$. Given any polarized K\"ahler metric $g$,
there is a Hermitian metric $h_L$ on $L$ whose Ricci form is equal to the K\"ahler form
$\omega_g$. 

The Bergman kernel, defined as a sequence of smooth (bundle-valued)   functions, plays  one of the central roles in K\"ahler-Einstein geometry.  In~\cite{t5}, Tian proved  the  $\mathcal C^2$-convergence theorem of the Bergman metric (whose K\"ahler potential is the Bergman kernel). 
A far-reaching generalization of Tian's  theorem was obtained by Catlin~\cite{cat} and Zelditch~\cite{zel}, using the existence of the parametrix of  Bergman kernel or Szeg\H o kernel (cf.~\cite{bs-1}):

\begin{theorem}[Catlin, Zelditch]\label{thm11}
Let $M$ be a compact complex manifold of dimension $n$ (over
$\mathbb{C}$) and let $(L,h_L)\rightarrow M$ be a positive Hermitian
holomorphic line bundle and $(E,h_E)$ a Hermitian vector bundle of
rank $r$. Let $g$ be the K\"ahler metric on $M$ corresponding to the
K\"ahler form $\omega_{g}=Ric(h_L)$.   Then there is an asymptotic
expansion of the Bergman kernel $\mathfrak B_m(x)$:
\begin{equation}\label{fud1}
\mathfrak B_m(x) \sim
a_0(x)m^n+a_1(x)m^{n-1}+a_2(x)m^{n-2}+\cdots
\end{equation}
for certain smooth coefficients $a_j(x)\in Hom(E,E)$ with $a_0=I$. More
precisely, for any $s$, the following inequalities hold:
\begin{equation}\label{main}
\|\mathfrak B_m(x) 
-\sum_{k=0}^sa_k(x)m^{n-k}\|_{\mathcal C^\mu}\leq C_{s,\mu}m^{n-s-1},
\end{equation}
where $C_{s,\mu}$ depends on $s,\mu$, the manifold $M$,  and the bundles $L,E$.
\end{theorem}

In~\cite{lu10} the first three coefficients were computed together with a general algorithm of computing any coefficients. 

There are several different proofs of the above expansion  theorem after Catlin and Zelditch. In~\cite{bbs}, Berman, Berndtsson and Sj\"ostrand gave a direct, constructive  approach, avoiding using the paramatrix of Bergman kernel. In~\cite{dlm}, Dai, Liu and Ma gave a heat kernel approach of the expansion.    In Englis~\cite{englis}, the first four coefficients were computed. More recently, Xu~\cite{haoxu} gave a graph-theoretic interpretation of the coefficients of the expansion. 

Because of the importance of the Bergman kernel expansion in complex geometry, we would like to give a sharper estimate of the   Bergman  expansion which is particularly useful in   the collapsing case~\cite{cliu-lu}.

In this paper, we  first prove an abstract version of the Bergman kernel expansion which is of its own interest. One of the innovations we obtain is that by systematically using the $K$-coordinate systems and $K$-frames, we are able to prove, in the abstract setting, that $\mathcal C^0$ expansion implies $\mathcal C^\mu$ expansion. Consequently, 
we generalize Theorem~\ref{thm11}.

\begin{theorem}\label{main-2}
We use the same notations as in Theorem~\ref{thm11}. Let $\eps>0$ be an absolute number (for example, we can take $\eps=1/16$).
Let $S_1,\cdots, S_k$ be peak sections, where $k=[\eps(\log m)^2]$.  Let $\mathfrak B_{m, peak}^k$  be the Bergman kernel with respect to  the peak sections   $S_1,\cdots, S_k$ (see Definition~\ref{def25} and Definition~\ref{def31}).  Then
\begin{enumerate}
\item $\mathfrak B_{m, peak}^k$ has an $\mathcal{C}^\infty$ asymptotic expansion;
\item The expansion stables to the Bergman kernel  expansion. That is,
\[
\|\mathfrak B_m(x)-\mathfrak B_{m, peak}^k(x)\|_{\mathcal C^\mu}\leq \frac{C}{m^{s-n+1}}
\]
for any $s\leq [\eps(\log m)^2]$.  
\end{enumerate}
In particular, for any $\mu\geq 0$, and any  increasing sequence  $\beta(m,\mu)\to\infty$, there exists an increasing integer sequence $\alpha(m,\mu)\to\infty$ such that 

\begin{equation}\label{s1.3}
\left\|\mathfrak B_m(x)-\sum_{k=0}^{\alpha(m,\mu)} a_k m^{n-k}\right\|_{\mathcal C^\mu}\leq \frac{\beta(m,\mu)}{m^{\alpha(m,\mu)-n+1}}.
\end{equation}

\end{theorem}

Since $\alpha(m,\mu)\to\infty$, for any fixed $s$, if $m$ is large enough, we must have $\alpha(m,\mu)>s$. Therefore, our theorem is a generalization of Theorem~\ref{thm11}.

The stretch from $(s+1)$ terms in Theorem~\ref{thm11} to $(\alpha(m,\mu)+1)$ terms in Theorem~\ref{main-2} provides a neat estimate in the real analytic case. 
When  the metrics of   $M, L, E$ are real analytic, 
  we are able to  improve~\eqref{s1.3}. More precisely, we prove the following result:

\begin{theorem}\label{thm1.2}
With the notations as in the above theorem, suppose that the
Hermitian metrics $h_{L}$ and $h_E$ are real analytic at a fixed
point $x$. Then the expansion
\[
\sum_{j=0}^\infty a_j(x)m^{-j}
\]
is convergent in $\mathcal C^\mu (\mu\geq 0)$ for $m$ large. Moreover,
we have
\begin{equation}\label{anal}
\|\mathfrak{B}_m(x)-\sum_{k=0}^\infty a_k(x) m^{n-k}\|_{\mathcal{C}^{\mu}}\leq
m^ne^{-\eps(\log m)^3}
\end{equation}
for some absolute constant $\eps>0$.
\end{theorem}

\begin{rem}
Both of  the above two theorems are shaper than Theorem~\ref{thm11}. But it  is not clear to the authors that  estimate ~\eqref{anal} is optimal. In~\cite{cliu1}, the first author proved the best possible estimate $e^{-\eps\, m}$ for compact Riemann surface (of genus $g\geq 2$) of constant Gauss curvature. But this is only a very special case.
\end{rem}

The following corollary of the above theorem provides bounds of the Bergman kernel coefficients, which is of independent interest.

\begin{cor}\label{cor11}  Assume that the metrics of $M,L,E$ are real analytic. Then there exists a constant $C>1$ such that the coefficients $a_j$ of the Bergman kernel expansion satisfy
\[
\|a_j(x)\|_{\mathcal C^\mu}\leq C^j
\]
for all $j\geq 0$.
\end{cor}

We remark that using the method in Lu-Shiffman (cf. ~\cite{lu-shiffman}*{Lemma 2.5}), similar results of this paper also hold for off-diagonal Bergman kernel.

 In~\cite{bbs}, a local reproducing kernel and then the local Bergman kernel were constructed explicitly. A relation between the local Bergman kernel and the Bergman kernel with respect to the peak sections (see Definiton~\ref{def25}), and its relation to the heat kernel paramatrix would provide insights of family version of the Catlin-Zelditch expansion. 
In particular, in both paper~\cite{bbs} and our paper, we only require the existence of a ``good'' holomorphic coordinate system, but we don't explicitly require that the injectivity radius to have  a lower bound. This observation may prove to be useful when we study the Bergman kernels on a family of \ka manifolds.

The method in proving the main results of this paper is  completely elementary.  It is similar to that in Shiffman~\cite{shiffman}, where we have to deal with matrices with increasing size.
We believe that our methods, including the discussion of the family version of $K$-coordinates and $K$-frames in the Appendix,  are  useful in many other places in complex geometry.

The paper is organized as follows. In Section~\S~\ref{s-2}, we prove Theorem~\ref{thm21}, an abstract version of the Bergman kernel  expansion.
The main results of this paper follow naturally from this abstract version and the hard analysis in \S ~\ref{3}-\ref{4}.
 In Section~\S~\ref{3}, we define the peak sections and verify assumption (1) of Theorem~\ref{thm21}; in Section~\S\ref{4}, we verify assumptions (2), (3), (4) of Theorem~\ref{thm21}; and in the Appendix, we include the discussion of family version of $K$-coordinates and $K$-frames.\\

{\bf Acknowledgement.} The first author thanks Chin-Lung Wang for his encouragement during the  preparation of this paper. She also thanks Chang-Shou Lin for his support during her stay in TIMS. The second author thanks Kengo Hirachi for his interest in our work and the discussions on the topic for many years.

\section{Preliminaries}
Throughout this paper, we use $C$ as a constant, which may differ from line to line. 
Let $(M,g)$ be a  K\"ahler manifold, and 
let $(E,h_E)$ be a Hermitian vector bundle of rank $r$.  For
$U,V\in\Gamma(M, E)$, the pointwise and the $L^2$
inner products are defined  as
\[
\langle U(x),V(x)\rangle_{h_{E}}
\] and
\numberwithin{equation}{section}
\begin{equation}
(U,V)=\int_{M}\langle
U(x),V(x)\rangle_{
h_E}dV_g,\label{1.1}
\end{equation} respectively, where
$dV_g={\omega_g^n}/{n!}$ is the volume form of $g$.

 Assume that $T_1,\cdots,T_d$ is an orthonormal basis of $H^0_{L^2}(M,E)$. Let  
\begin{equation}\label{97}
(e_1,\cdots,e_r)
\end{equation}
 be holomorphic local frames of $E$. 
Let 
\[
T_j=\sum_{\alpha=1}^r a_{j\alpha}e_\alpha
\]
for $1\leq j\leq d$.
Then $a_{j\alpha}$ are local holomorphic  functions. Let $A(x)$ denote the $d\times r$ matrix $(a_{j\alpha}(x))$. The Bergman kernel,  which   is an element of ${\rm Hom}\, (E_y,E_x)$  for any $x,y\in M$,  is defined by 
\[
\mathfrak B(x,y)=H(y)A^*(y)A(x),
\]
where $H=(h_{\alpha\bar\beta})=(\langle e_\alpha,e_\beta\rangle)$ is the metric matrix and $A^*$ is the complex conjugate of $A$.

 The name {\it Bergman kernel} is justified by the following property. For any $f\in \Gamma_{L^2}(M,E)$, if we write $f$ as a row vector under the frame~\eqref{97}, then 
 \[
 \int_M f(y)\mathfrak B(x,y) dV_y
 \]
 is the orthogonal projection  from $\Gamma_{L^2}(M,E)$ to   $H^0(M,E)$.

In this paper, we study $\mathfrak B(x)=\mathfrak B (x,x)=HA^*A$, which we  also call the Bergman kernel.

The Bergman kernel  $\mathfrak B(x)$ is independent of the choice of orthonormal basis $T_1,\cdots, T_d$. 
Let $S_1,\cdots,S_d$ be {\it any} basis of $H^0(M,E)$. Let 
\begin{equation}\label{ff}
F=(F_{ij})=(S_i,S_j).
\end{equation}
Let $P$ be a matrix such that $PFP^\ast=I$. 
If we write
\begin{equation}\label{p0}
S_j=\sum_{\alpha=1}^rb_{j\alpha}e_\alpha,
\end{equation}
then the Bergman kernel can be represented by
\begin{equation}\label{berg-3}
\mathfrak{B}(x)=H(PB)^\ast PB=HB^\ast F^{-1}B,
\end{equation}
where $B=(b_{j\alpha})$ is a local $d\times r$ matrix-valued function.
\\

Now we assume that the K\"ahler manifold $M$ is {\it polarized}. That is, there is an ample Hermitian line bundle $(L,h_L)$ over $M$ such that the K\"ahler metric $g$ is the curvature ${\rm Ric}(h_L)$. 
Let $m$ be a large positive integer. Let $\mathfrak B_m(x)$ be the Bergman kernel of the bundle $L^m\otimes E$, which is 
a sequence of $Hom(L^m\otimes E,L^m\otimes E)(=Hom(E,E))$-valued smooth functions. The main purpose of this paper is to study the asymptotic expansion of these functions. 

\section{An abstract version of Bergman kernel expansion}\label{s-2}
We begin with some abstract  discussions of the bundle-valued function expansions. 
\begin{definition}\label{def2.1}
We say a sequence of matrix-valued functions $f_m(x)$ has a strongly $\mathcal C^\mu$   asymptotic expansion, if there exists a sequence $\sigma(m)\to\infty$ for $m\to\infty$, such that for any $s<\sigma(m)$ and $\mu\geq 0$,  we have
\begin{equation}\label{1.9.1}
\left\|f_m(x)-m^n\left(a_0(x)+
\frac{a_1(x)}{m}+\cdots+\frac{a_s(x)}{m^s}\right)\right\|_{\mathcal C^\mu}\leq \frac{C(s,\mu)}{m^{s-n+1}}
\end{equation}
for  matrix-valued functions 
$a_0(x)$, $\cdots$, $a_s(x)$, $\cdots$, and a constant  $C(s,\mu)$ independent to $m$. If the metrics of $M,L,E$ are real analytic, then we further assume that there exists a constant $C_1>1$, independent to $s$ and $x$, such that $C(s,\mu)\leq C_1^{s+1}$.
\end{definition}

It is natural to believe that   a sequence of  functions $f_m(x)$ has a strongly $\mathcal C^\mu$   asymptotic expansion if and only if at any point of $M$, it has a strongly $\mathcal C^\mu$   asymptotic expansion. However, we {\it can't} prove such a statement. 
Since the $\mathcal C^\mu$-norm depends on the choices of local coordinates and frames, in order to obtain~\eqref{1.9.1}, we need to prove the existences of a smooth family of  ``good'' coordinates and frames at any point of the manifold. We first  make the following definition.

\begin{definition}\label{def32-2}
Let $\mu$ be a nonnegative  integer.
We say a sequence of matrix-valued functions $f_m(x)$ has a strongly  $\mathcal C^\mu$ asymptotic expansion at  any point, 
 if there exists  a constant $C$ independent to $m$ such that for  any  $x_0\in M$,
 there exist local $K$-coordinate system and $K$-frames at $x_0$ of order $p>\mu+\sigma(m)$,
and matrix-valued functions\footnote{Strictly speaking, these functions are vector-valued functions whose components are matrix-valued functions. This is because a $\mu$-th derivative has many components. But for the sake of simplicity, we don't emphasis this fact.}  $a^\mu_0(x)$, $\cdots$, $a^\mu_s(x)$, $\cdots$ in a neighborhood of $x_0$ with the following property: for any  $s<\sigma(m)$, $\mu'\leq \mu$, there exist constants $C(s,\mu')$ which are independent to $m$ and $s$ such that
 \begin{equation}\label{1.9}
\left|D_{x_0}^{\mu'}(f_m(x))-m^n\left(a^{\mu'}_0(x)+
\frac{a^{\mu'}_1(x)}{m}+\cdots+\frac{a^{\mu'}_s(x)}{m^s}\right)\right|(x_0)\leq \frac{C(s,{\mu'})}{m^{s-n+1}},
\end{equation}
where  $D_{x_0}^{\mu'}$ is the ${\mu'}$-th derivative with respect to the above-mentioned $K$-coordinate system and frames centered at $x_0$.

In particular, we say $f_m(x)$ has a strongly  $\mathcal{C}^0$ asymptotic expansion\footnote{By this notation, we have $a_j^0=a_j$ for any $j\geq 0$.} at any point, if  for any $s<\sigma(m)$
\begin{equation}\label{1.9.2}
\left|f_m(x_0)-m^n\left(a_0(x_0)+
\frac{a_1(x_0)}{m}+\cdots+\frac{a_s(x_0)}{m^s}\right)\right|\leq \frac{C(s)}{m^{s-n+1}}
\end{equation}
for any $x_0\in M$,
where $C(s)$ are constants  independent to $m$ and $x_0$.

Similar to the previous definition,  if the metrics of $M,L,E$ are real analytic, then we further assume that there exists a constant $C_1>1$, independent to $s$, such that $C\leq C_1^{s+1}$.

\end{definition}

\begin{lemma}\label{lem2.1}
For any  positive  integer $\mu$, a sequence of matrix-valued functions $f_m$ has  a  strongly $\mathcal C^{\mu}$ asymptotic expansion,  if  $f_m$ has a   strongly $\mathcal C^\mu$ asymptotic expansion at  any point.
\end{lemma}

\begin{proof}
We claim that for any nonnegative integers $\mu,s$, we have\footnote{It should be noted that these functions are not even continuous if we choose a discontinuous family of $K$-coordinates.}
\[
a_s^\mu(x)=D_x^\mu a_s(x),
\]
where $D_x^{\mu'}$ is the ${\mu'}$-th derivative with respect to the above-mentioned $K$-coordinate system and frames centered at $x$.

For $\mu=0$, the claim follows from definition. Inductively, assume that the claim is proved to be true for any $\mu'\leq \mu-1$ and 
for any $s'\leq s$ when $\mu'=\mu$. We shall prove that 
\begin{equation}\label{claim-1}
a_{s+1}^\mu(x)=D_x^\mu a_{s+1}(x).
\end{equation}
Without loss of generality, we only need to prove the above statement on an open neighborhood $U$ of a fixed point $x_0$. By  Corollary~\ref{corA1}, and by shrinking $U$ if necessary,  there exists  a holomorphic family of $K$-coordinate systems on $U$. In what follows, we shall fix such a family of coordinate systems.

We define a sequence of functions
\[
\xi_{m,s}(x)=m^{s+1}\left(\frac{f_m(x)}{m^n}-\left(a_0(x)+\cdots+\frac{a_s(x)}{m^s}\right)\right).
\]
By ~\eqref{1.9} and the inductive assumption, 
$D^\mu_x(\xi_{m,s}(x))$
are uniformly convergent to $a_{s+1}^\mu(x)$ on $U$.
 Using  the chain rule,  for any $\mu'\leq\mu$, we have
\begin{equation}\label{chain-rule}
D_{x_0}^{\mu'}(\xi_{m,s}(x))=F_{\mu'}(D_{x}^{\mu'}(\xi_{m,s}(x)), D_{x}^{\mu'-1}(\xi_{m,s}(x)),\cdots,
\xi_{m,s}(x))
\end{equation}
for  functions $F_{\mu'}$ depending only on the family of $K$-coordinate systems. Therefore, the functions $\{\xi_{m,s}(x)\}$ are convergent on $U$ in the $\mathcal C^\mu$-norm to some function $\tilde a_{s+1}(x)$. By taking limit in
~\eqref{chain-rule} and the inductive assumption, we have
\[
D_{x_0}^{\mu}(\tilde a_{s+1}(x))=F_{\mu}(a^\mu_{s+1}(x), D_{x}^{\mu-1}(a_{s+1}(x)),\cdots,
a_{s+1}(x)).
\]
On the other hand, by the chain rule again, we have
\[
D_{x_0}^{\mu}(\tilde a_{s+1}(x))=F_{\mu}(D_{x}^{\mu}(a_{s+1}(x)), D_{x}^{\mu-1}(a_{s+1}(x)),\cdots,
a_{s+1}(x)).
\]
By further shrinking $U$ if necessary, we can represent both $a^\mu_{s+1}(x)$ and  $D_{x}^{\mu}(a_{s+1}(x))$ as the exact same polynomial of 
\[
D_{x_0}^{\mu}(\tilde a_{s+1}(x)), D_{x}^{\mu-1}(a_{s+1}(x)),\cdots,
a_{s+1}(x).
\]
Therefore they must be equal
\[
a^\mu_{s+1}=D_{x}^{\mu}(a_{s+1}(x)),
\]
and the claim is proved. 

Lemma~\ref{lem2.1} then follows from~\eqref{1.9} and ~\eqref{chain-rule}.

\end{proof}

\label{page-6}
Let $P=(p_1,\cdots,p_n)$ be a multiple index and let $1\leq \alpha\leq r$. Define the lexicographical order on the set of $(P,\alpha)$'s. That is, $(P,\alpha)<(Q,\beta)$ if 
\begin{enumerate}
\item  $\sum p_i<\sum q_i$, or
\item  $p_1=q_1,\cdots,p_\ell=q_\ell$ but $p_{\ell+1}<q_{\ell+1}$ for some $0\leq \ell\leq n-1$, or
\item  $P=Q$, but $\alpha<\beta$.
\end{enumerate}
Such an order gives rise to the function $j=\sigma(P,\alpha)$. For example, $1=\sigma((0,\cdots,0),1)$, $2r+2=\sigma((0,1,\cdots, 0),2)$, etc. Its inverse function is defined by $(P,\alpha)=(R(j),\zeta(j))$. For an example, $(R(1),\zeta(1))=((0,\cdots,0),1)$ and $(R(2r+2),\zeta(2r+2))=((0,1,\cdots,0),2)$, etc.

\begin{definition}\label{def23-1}
Let $\sigma(m)$ be  sequences of  positive integers such that $\sigma(m)\to\infty$ as $m\to\infty$. Let $x_0\in M$ be a fixed point. At $x_0$, let $(z_1,\cdots,z_n)$ be  $K$-coordinates  and  $e_L$, $(e_1,\cdots, e_r)$ be $K$-frames for $L$, $E$, respectively (of order $\mu+\sigma(m)$ if the metrics are smooth and  of order $+\infty$ if the metrics are analytic). 
Let $S_1,\cdots,  S_d$ be a basis of $H^0(M,L^m\otimes E)$. We say  that it is  a  strongly regular basis at $x_0$,  if for the smallest $k$ such that $|R(k)|\geq\sigma(m)$, \begin{enumerate}
\item for $1\leq j\leq k$, $S_j(z)=z^{R(j)}e_L^m\otimes e_{\zeta(j)}+o(|z|^{\sigma(m)})$;
\item for $j>k$, $S_j(z)=o(|z|^{\sigma(m)})$.
\end{enumerate}
\end{definition}

\begin{rem}
The motivation of the above definition is the existence of peak sections (Definition~\ref{def31}).
\end{rem}

For the rest of the paper, we shall let
\[
\sigma(m)=s=[\eps(\log m)^2],
\]
where $\eps>0$ is a  small absolute positive number. 
Let $\delta_{-1}=0$ and 
\[
\delta_j=r\binom{n-1+j}{n-1}
\]
 for $j=0,\cdots,s$. Then $\delta_j$ is the number of indices $(P,\alpha)$ such that $|P|=j$ and $1\leq\alpha \leq r$.
Let\footnote{By math induction, $\delta_{s+1}=r\binom{n+s}{n}$.}
\begin{align*}
&\delta_{s+1}=r\sum_{j=1}^{s}\delta_j.\\
\end{align*}

Then we have 
\begin{equation}\label{plk-1}
k=r\sum_{j=1}^{s}\delta_j \leq 2rs^n.
\end{equation}
In particular,
\[
|R(k)|=s.
\]

\begin{definition}\label{24-2}
 A strongly regular basis $\{S_j\}_{(j=1,\cdots d)}$ is called almost orthonormal, if
\begin{equation}\label{nmb}
\frac{|(S_i,S_j)|}{\|S_i\|\cdot\|S_j\|}\leq \frac{C}{m^{1+(|R(j)|-|R(i)|)/2}}
\end{equation}
for any $i< j\leq 2k$, and 
\[
\begin{array}{ll}
(S_i,S_j)=0& i\leq k, j>2k;
\\
(S_i,S_j)=\delta_{ij}& i,j>2k,
\end{array}
\]
where the constant $C$ is independent to $m$ and $x_0$.  Moreover, if the metrics of $M,L,E$ are real  analytic, then we assume that there exists a constant $C_1>1$, independent to $m,x_0$ and $i,j$, such that
\[
C\leq C_1^{|R(j)|-|R(i)|}.
\]
\end{definition}

 Let $F=((S_i,S_j))_{i,j=1,\cdots d}$ be the $d\times d$  metric matrix. 
For an almost orthonormal  basis, we can write $F$ as
 \[
 F=\begin{pmatrix}
F_1\\
&I_{d-2k}
 \end{pmatrix},
 \]
 where $F_1$ is a $(2k)\times (2k)$
 matrix. In order to study the inverse matrix of $F$, we only need to study the inverse matrix of $F_1$.
 
 In what follows, we shall prove our key observation: for an almost orthonormal basis, up to order $[\eps(\log m)^2]$,  the  expansions of the $k\times k$ minor of  $F^{-1}$ only depends on $(S_i,S_j)$ for $i,j\leq k$ and are independent to the rest of the matrix!  Since we shall pick $(S_1,\cdots, S_k)$ as {peak} sections  later,  we reduce the expansion of Catlin-Zelditch ($\mathfrak{B}_m$) into the expansions $\mathfrak{B}_{m, peak}^k$. That is, we reduce the problem to find the expansion of the inverse of the inner products of the peak sections.
 
 Let 
 \[
 F_2=\Lambda^{-1}F_1\Lambda^{-1},
 \]
 where 
 \[
 \Lambda=\begin{pmatrix}
 \|S_1\|\\
 &\ddots\\
 &&\|S_k\|\\
 &&&1\\
 &&&&\ddots\\
 &&&&&1
 \end{pmatrix}.
 \]

 We  represent the matrix $F_2$ as a block matrix.  For  any $0\leq\xi,\eta\leq s+1$, define the matrix $A_{\xi\eta}$ to be the $\delta_\xi\times\delta_\eta$ matrix whose entries are\footnote{By this definition, $A_{s+1,s+2}=0$ and $A_{s+2,s+1}=0$.}
 \[
 \delta_{ij}-\frac{(S_i,S_j)}{\|S_i\|\cdot \|S_j\|}
 \]
 for 
 \[
 \sum_{\gamma=-1}^{\xi-1}\delta_\gamma<i\leq   \sum_{\gamma=-1}^{\xi}\delta_\gamma;
 \quad
  \sum_{\gamma=-1}^{\eta-1}\delta_\gamma<j \leq   \sum_{\gamma=-1}^{\eta}\delta_\gamma.
 \]
 Let $A$ be the block matrix defined by
 \[
 A=(A_{\xi\eta})_{(0\leq\xi,\eta\leq s+1)}.
 \]
 Then we have
 \[
 F_2=I-A.
 \]

 Define $\|A\|$ to be the maximum entries  of the matrix $A$. Let ${\rm column}\,(A)$ be the number of the columns of the matrix $A$. Then we have
 \begin{align}\label{norm-3}
 \begin{split}
& \|A+B\|\leq \|A\|+\|B\|;\\
 &
 \|AB\|\leq {\rm column}\,(A)\cdot \|A\|\cdot\|B\|.
 \end{split}
 \end{align}

By~\eqref{nmb}, we have
\begin{equation}\label{norm-4}
\|A_{\xi\eta}\|\leq\frac{C}{m^{1+\frac 12|\xi-\eta|}}.
\end{equation}
Using the expansion
\[
F_2^{-1}=I+A+A^2+\cdots,
\]
for any fixed $(\xi_0,\eta_0)$, we have

\[
(F_2^{-1}-I)_{\xi_0\eta_0}=\sum_{k=1}^\infty\sum_{i_1,\cdots,i_{k-1}}A_{\xi_0 i_1}A_{i_1i_2}\cdots A_{i_{k-1} \eta_0}.
\]

By ~\eqref{plk-1}, ${\rm column}\,(A)\leq 2rs^n$. Therefore by~\eqref{plk-1} and ~\eqref{norm-3}, for fixed $\xi_0,i_1,\cdots, i_{k-1},\eta_0$, we have
\[
\left\|A_{\xi_0 i_1}A_{i_1i_2}\cdots A_{i_{k-1}\eta_0}\right\|\leq \left(\frac{C\cdot 4 rs^n}{m}\right)^k.
\]
Therefore
\[
\left\|\sum_{i_1,\cdots,i_{k-1}}A_{\xi_0 i_1}A_{\xi_0 i_1}A_{i_1i_2}\cdots A_{i_{k-1}\eta_0}\right\|\leq (s+1)^k\left(\frac{C\cdot 4rs^n}{m}\right)^k.
\]
For fixed $s$, if $m$ is large enough, then we have
\[
\left\|\sum_{k=s+1}^\infty\sum_{i_1,\cdots,i_{k-1}}A_{\xi_0 i_1}A_{i_1i_2}\cdots A_{i_{k-1}\eta_0}\right\|\leq \frac{C}{m^{s+1}}.
\]
 Similarly, we consider the  terms
\[
\sum_{k=1}^s\sum_{\text{some } i_j=s+1}A_{\xi_0 i_1}A_{i_1i_2}\cdots A_{i_{k-1}\eta_0}.
\]
If some $i_j=s+1$, we must have
\[
|\xi_0-i_1|+|i_1-i_2|+\cdots+|i_{k-1}-\eta_0|\geq 2s+2-\xi_0-\eta_0.
\]
Thus by~\eqref{norm-4} we have
\[
\left\|\sum_{k=1}^s\sum_{\text{some } i_j=s+1}A_{\xi_0 i_1}A_{i_1i_2}\cdots A_{i_{k-1}\eta_0}\right\|\leq \frac{C}{m^{s+1-\frac 12(\xi_0+\eta_0)}}.
\]
Consequently, we have
\begin{align}\label{5.2}
\begin{split}
&\left\|(F_2^{-1})_{\xi_0\eta_0}-\left(I_{\xi_0\eta_0}+\sum_{k=1}^s\sum_{i_1,\cdots,i_{k-1}=0}^s A_{\xi_0 i_1}A_{i_1i_2}\cdots A_{i_{k-1}\eta_0}\right)\right\|\\ &\leq \frac{C}{m^{s+1-\frac{1}{2}(\xi_0+\eta_0)}}.
\end{split}
\end{align}
If the metrics of $M,L,E$ are real analytic, then using the same estimate as above, we will get\footnote{In the analytic case, $s\to\infty$ as $m\to\infty$.}
\begin{align}\label{5.2-2}
\begin{split}
&\left\|(F_2^{-1})_{\xi_0\eta_0}-\left(I_{\xi_0\eta_0}+\sum_{k=1}^s\sum_{i_1,\cdots,i_{k-1}=0}^s A_{\xi_0 i_1}A_{i_1i_2}\cdots A_{i_{k-1}\eta_0}\right)\right\|\\ &\leq \left(\frac{C_1(\log m)^{4n}}{m}\right)^{s+1-\frac{1}{2}(\xi_0+\eta_0)}
\end{split}
\end{align}
for $s=[\eps(\log m)^2]$.

 Let 
\[
\delta(\xi_0,\eta_0,i_1,\cdots, i_{k-1})=|\xi_0-i_1|+|i_1-i_2|+\cdots +|i_{k-1}-\eta_0|.
\]
We use the notation
$f(m)\ll g(m)$ to denote that 
\begin{enumerate}
\item both $f(m), g(m)$ are Taylor  series of $m^{-1}$;
\item the coefficients of $g(m)$ are non-negative;
\item the coefficients of $g(m)-f(m)$ are non-negative. 
\end{enumerate}
By the assumption of $A_{ij}$, we have
\[
\|A_{ij}\|\ll \left(\frac Cm\right)^{1+\frac{|i-j|}{2}}\cdot\frac{1}{1-\frac Cm}
\]
for some constant $C>1$.
It follows that
\begin{align}\label{34-2}
\begin{split}
&\sum_{k=1}^s\sum_{i_1,\cdots,i_{k-1}=0}^s A_{\xi_0 i_1}A_{i_1i_2}\cdots A_{i_{k-1}\eta_0}\\
&\ll\sum_{k=1}^s \sum_{i_1,\cdots, i_{k-1}=0}^s \left(\frac Cm\right)^{k+\delta(\xi_0,\eta_0,i_1,\cdots, i_{k-1})}\left(1-\frac Cm\right)^{-k}.
\end{split}
\end{align}

Let $\#\delta$ be the number of $i_1,\cdots, i_{k-1}$ such that $\delta(\xi_0,\eta_0,i_1,\cdots, i_{k-1})=\delta$. In the expansion of the left hand side of ~\eqref{34-2}, the coefficient of $m^{-b}$ is not more than
\[
C^b\sum_{\delta+k\leq b}\frac{(b-\delta)!}{k!(b-\delta-k)!}\cdot\#\delta.
\]
By Lemma~\ref{abcd-2} below, the coefficients of $m^{-b}$ is not more than
\[
C^b\sum_{\delta+k\leq b} \frac{(b-\delta)!}{k!(b-\delta-k)!}\cdot\#\delta\leqq C^b
\]
for a possibly larger  constant $C$.

 Replacing $s$ in ~\eqref{5.2-2} by $s+1+\frac 12(\xi_0+\eta_0)$, we get
 \begin{align}\label{5.2-3}
 \begin{split}
 &\left\|(F_2^{-1})_{\xi_0\eta_0}-\left(I_{\xi_0\eta_0}+\sum_{k=1}^{s+1+\frac 12(\xi_0+\eta_0)}\sum_{i_1,\cdots,i_{k-1}=0}^{s+1+\frac 12(\xi_0+\eta_0)} A_{\xi_0 i_1}A_{i_1i_2}\cdots A_{i_{k-1}\eta_0}\right)\right\|\\ &\leq \left(\frac{C_1(\log m)^{4n}}{m}\right)^{s+2}.
 \end{split}
 \end{align}
 
 By the estimate on the coefficients of $m^{-b}$, we obtain
 \begin{align}\label{5.2-4}
 \begin{split}
& \|(F_2^{-1})_{\xi_0\eta_0}-\text{ the expansion up to order $s$} \|\\
&\leq \left(\frac Cm\right)^{s+1}+
 \left(\frac{C_1(\log m)^{4n}}{m}\right)^{s+2}\leq \left(\frac Cm\right)^{s+1}.
 \end{split}
 \end{align}

By using the same method, when the metrics are smooth, from~\eqref{5.2}, we obtain
\begin{equation}\label{5.2-5}
\|(F_2^{-1})_{\xi_0\eta_0}-\text{ the expansion up to order $s$} \|\leq\frac{C}{m^{s+1}},
\end{equation}
 where $C$ is a constant independent to $m$ (but may depend on $s$). \\

 We need the following combinatorial  lemma in completing the above estimate.

\begin{lemma}\label{abcd-2}
There exists a constant $C>1$ such that
\[
\#\delta\leq C^{\delta+k}.
\]
\end{lemma}

\begin{proof} It is well-known that the number of non-negative solutions of the equation
\[
a_1+\cdots+a_k=\delta
\]
is equal to 
\[
\frac{(\delta+k-1)!}{(k-1)!\delta!}\leq C^{\delta+k}.
\]
On the other hand, $\#\delta$ is not more than  the number of solutions of the equations $i_{j-1}-i_{j}=\pm a_{j}$ for $1\leq j\leq k$ (where we define $i_0=\xi_0, i_{k}=\eta_0$). Thus we have the estimate
\[
\#\delta\leq 2^k C^{\delta+k}\leq C^{\delta+k}
\]
for some larger constant $C>1$.

\end{proof}

Now we state the  following abstract version of the Bergman kernel expansion:

\begin{theorem}\label{thm21} Let $\sigma(m)=[\eps(\log m)^2]$ for a sufficiently small absolute number  $\eps>0$.  Let $k$ be the smallest  positive integer such that $|R(k)|\geq\sigma(m)$. Let $s=|R(k)|$.
Assume that for any $x_0\in M$, there exists a basis $\{S_j\}_{(j=1,\cdots,d)}$ of $H^0(M,E\otimes L^m)$ such that 
\begin{enumerate}
\item it  is a strongly regular basis in the sense of Definition~\ref{def23-1};
\item it  is almost orthonormal in the sense of Definition~\ref{24-2};
\item the expressions
\[
\frac{(S_i,S_j)}{\|S_i\|\cdot\|S_j\|}
\]
have strongly $\mathcal C^0$ asymptotic expansions (in the sense of ~\eqref{1.9.2})
\[
\frac{(S_i,S_j)}{\|S_i\|\cdot\|S_j\|}\sim \frac{r_l}{m^l}+\frac{r_{l+1}}{m^{l+1}}+\cdots
\]
for $|R(i)|+|R(j)|\leq \mu$ with 
 $l\geq 1+(|R(j)|-|R(i)|)/2$.
Moreover, we have $|r_{t}|\leq C$, where $C$ is independent to $i,j$ for all $l\leq t\leq s$.  If the metrics of $M,L,E$ are real analytic, we further
assume that there exists a constant $C_1>1$, independent to $m,x_0$ and $i,j$, such that $C\leq C_1^t$.
\item the expressions
$\|S_i\|$
have strongly $\mathcal C^0$  asymptotic expansions  (in the sense of ~\eqref{1.9.2})
\[
\|S_i\|\sim \frac{r_l}{m^l}+\frac{r_{l+1}}{m^{l+1}}+\cdots
\]
for $|R(i)|\leq \mu$, where 
 $l=(n+|R(i)|)/2$, and $\|r_l\|_{\mathcal C^\mu}\geq c>0$.
Moreover, we have $|r_{t}|\leq C$, where $C$ is independent to $i,j$ for all $l\leq t\leq s$.  If the metrics of $M,L,E$ are real analytic, we further
assume that there exists a constant $C_1>1$, independent to $m,x_0$ and $i,j$, such that $C\leq C_1^t$.
\end{enumerate}
 Then we have the following  Catlin-Zelditch type expansion:
 For any increasing sequence $\beta(m,\mu)\to\infty$, there exists a sequence $\alpha(m,\mu)\to\infty$ such that
\begin{equation}\label{real-2-0}
\left\|\mathfrak B_m-\sum_{k=0}^{\alpha(m,\mu)} a_k m^{n-k}\right\|_{\mathcal C^\mu}\leq \frac{\beta(m,\mu)}{m^{\alpha(m,\mu)-n+1}}.
\end{equation}
 Moreover, 
 if the metrics of $M,L,E$ are real analytic, then we have the following sharper estimate
\begin{equation}\label{real-2}
\left\|\mathfrak B_m-\sum_{k=0}^\infty a_k m^{n-k}\right\|_{\mathcal C^\mu}\leq m^ne^{-\eps(\log m)^3}
\end{equation}
for a small absolute constant $\eps$.
\end{theorem}

\begin{proof} Let $a$ and $H$ be the Hermitian metrics of $L$ and $E$, respectively.  Then by~\eqref{berg-3}, near a fixed point $x_0$,
\[
\mathfrak B_m=a^m HB^*F^{-1}B.
\]
It is apparent that any derivatives of the functions $a^m$ and $H$ at $x_0$ are  polynomials of $m$.
For the rest of the paper, we shall repeatedly use the following elementary fact:  if functions $f_m, g_m$ have asymptotic expansions, so is their product. 

Thus in order to prove the result, we only need to prove the existence of the strongly $\mathcal C^\mu$ asymptotic expansion of $B^*F^{-1}B$.
 
Using~\eqref{p0} and the definition of $S_j$, we have
\[
b_{j\alpha}=\delta_{\alpha,\zeta(j)}\, z^{R(j)}+o(|z|^{\sigma(m)}),
\]
where $\delta_{ij}$ is the Kronecker symbol.
It follows that 
\[
\left.\frac{\pa^P b_{j\alpha}}{\pa z^P}\right|_{z=0}=\delta_{\alpha,\zeta(j)}\,\delta_{P,R(j)} \, R(j)!
\]
for any $|P|\leq \mu<\sigma(m)$.

 By a straightforward computation, the $(i, j)$-th entry of the derivatives of the matrix $B^*F^{-1}B$ is
\begin{align}\label{xyz-1}
\begin{split}
&\left(\left.\frac{\pa^{|P|+|Q|} }{\pa z^P\pa\bar z^Q}(B^*F^{-1}B)\right|_{z=0}\right)_{\alpha\beta}=P!Q!(F^{-1})_{i j}\\&
=P!Q!(F_2^{-1})_{ij}\|S_i\|^{-1}\cdot\|S_j\|^{-1},\end{split}\end{align}
for any $1\leq \alpha,\beta\leq r$; $|P|+|Q|\leq\mu$,
where $i=\sigma(P,\alpha), j=\sigma(Q,\beta)$.  

By~\eqref{5.2-5} and assumption (4) in Theorem~\ref{thm21}, the above expression has a strongly $\mathcal C^0$ asymptotic expansion. This implies that the Bergman kernel has a $\mathcal C^\mu$ expansion at a point in the sense of Definition~\ref{def32-2}. By  Lemma~\ref{lem2.1}, for any $s<\sigma(m)-\mu$, ~\eqref{main} is valid. Without loss of generality, we assume that the constants $C_{s,\mu}$ in~\eqref{main} are increasing with respect to $s$ and is divergent to $\infty$ as $s\to\infty$.
Let $\beta(m,\mu)<\sigma(m)-\mu$ be an increasing positive sequence that is divergent to $\infty$ as $m\to\infty$. Let $\alpha(m,\mu)\to\infty$ be an  increasing sequence of integers such that 
\[
C_{\alpha(m,\mu),\mu}<\beta(m,\mu).
\]
 Then ~\eqref{real-2-0} is valid.

 If the metrics are real analytic,   we use~\eqref{5.2-4} for $s$ and $s-1$, we have
 \begin{align*}
 & \|(F_2^{-1})_{\xi_0\eta_0}-\text{ the expansion up to order $s$} \|\leq\left(\frac Cm\right)^{s+1};\\
 & \|(F_2^{-1})_{\xi_0\eta_0}-\text{ the expansion up to order $s-1$} \|\leq\left(\frac Cm\right)^{s}.
 \end{align*}
 Let $b_j$ be the coefficient of the $m^{-j}$ in the expansion of $F_2^{-1}$. From the above inequality, we have $\|b_j\|_{\mathcal C^\mu}\leq C^j$. Therefore, the coefficients $a_j$ of the Bergman kernel expansion satisfy
 \begin{equation}\label{cor11-1}
 \|a_j\|_{\mathcal C^\mu}\leq C^j
 \end{equation}
 for $j\geq 0$. 
Consequently, we have
 \[
 \|\sum_{k=s+1}^\infty a_k m^{n-k}\|_{\mathcal C^\mu}\leq m^n\left(\frac Cm\right)^{s+1}.
 \]
 On the other hand, since we already obtain the estimate
 
 \[
 \|\mathfrak B_m
-\sum_{k=0}^sa_jm^{n-k}\|_{\mathcal C^\mu}\leq m^n\left(\frac Cm\right)^{s+1},
\]
combining the above two  inequalities and  taking $s=[\eps(\log m)^2]$, we obtain~\eqref{real-2}.
 
\end{proof}

We end this section by making the following relative version of the Bergman kernel.

\begin{definition}\label{def25}
Let $\{S_1,\cdots, S_d\}$ be a basis of $H^0(M,L^m\otimes E)$. Then the Bergman kernel 
$\mathfrak B_{m, peak}^k$ with respect to part of the basis $\{S_1,\cdots, S_k\}$ is defined as
\[
\mathfrak B_{m, peak}^k=a^mH(B_k)^*F_k^{-1}B_k,
\]
where $B_k, F_k$ are defined similarly  to ~\eqref{ff} and~\eqref{p0}. That is, $F_k$ is the $k\times k$ matrix defined by
\[
(F_k)_{ij}=(S_i,S_j),
\]
and $B_k$ is the $k\times r$ matrix whose entries  $(b_k)_{j\alpha}$ are defined by 
\[
S_j=\sum_{\alpha=1}^r (b_k)_{j\alpha} e_\alpha
\]
for some basis $e_1,\cdots, e_r$.
\end{definition}

\section{Peak sections}\label{3}

The peak sections were introduced 
in~\cite{t5} (see also~\cite{lu10}). The results of this section are mostly known, except we make efforts to extend our estimate up to order $[\eps(\log m)^2]$, which grows to infinity as $m\gg 0$.

We begin with citing the following well-known result (cf. ~\cites{de, t5}). 
\begin{prop}\label{prop1}
Suppose that $(M,g)$ is an $n$-dimensional compact  K\"ahler
manifold,  $(L,h_L)$ is a Hermitian
line bundle over $M$, and
$(E,h_E)$ is a Hermitian vector bundle over $M$. We assume that ${\rm Ric}\,(h_L)=\omega_g$ defines the \ka metric of $M$. Let $\Theta_{E}$ be the curvature of
$E$ with respect to $h_{E}$ and $I_{E}$ be the identity map of $E$.
Let $\Psi$ be a function on $M$ which can be approximated by a
decreasing sequence of smooth function
$\{\Psi_{\ell}\}_{\ell=1}^{\infty}$. Assume that the endomorphisms
\[\Theta_{\ell}=\Theta_{E}+I_{E}\otimes\partial\bar{\partial}
\Psi_{\ell}+\frac{2\pi}{\sqrt{-1}}I_{E}\otimes (Ric(h)+Ric(g))\] on
$\bigwedge^{0,1}(M,L\otimes E)$ satisfy
\begin{equation}\label{psi}
\langle\Theta_{\ell}(\varphi),\varphi\rangle_{g^{\ast}\otimes
h_{L}\otimes h_{E}}\geq C \|\varphi\|^{2}_{g^{\ast}\otimes
h_{L}\otimes h_{E}}\end{equation} for $\varphi\in
\bigwedge^{0,1}(M,L\otimes E)$, where $C>0$ is a constant independent
of $\ell$. Then for any $w\in \bigwedge^{0,1}(M,L\otimes E)$ with
$\bar{\partial}w=0$ and
\[\int_{M}\|w\|^{2}_{g^{\ast}\otimes
h_{L}\otimes h_{E}}e^{-\Psi}dV_{g}< \infty,
\]
there exists $u\in
C^\infty(M,L\otimes E)$ such that $\bar{\partial}u=w$ and
\begin{equation}
\int_{M}\|u\|^{2}_{h_{L}\otimes h_{E}}e^{-\Psi}dV_{g}\leq
\frac{1}{C}\int_{M}\|w\|^{2}_{g^{\ast}\otimes h_{L}\otimes
h_{E}}e^{-\Psi}dV_{g}.\label{2.4}
\end{equation}
\end{prop}

\qed

\begin{rem}\label{rem21}
The choice of $u$ is, of course, not unique. In order to
construct the peak sections, we need to  fix the  solution of the equation
$\bar\pa u=w$. Let $\Delta$ be the $(0,1)$-Laplacian on  bundle
$L\otimes E$ with respect to the corresponding
metrics and the weight function $\Psi$. By~\eqref{psi} and the
Weitzenb\"ock formula, the smallest eigenvalue $\lambda_1$ of
$\Delta$ satisfies $\lambda_1\geq C$. It follows that the Green's
operator $G$ of $\Delta$ exists and is a bounded operator on $L^2$ spaces.

Let
\[
u=\bar\pa^* Gw.
\]
Then we have
\[
\bar\pa u=\bar\pa\bar\pa^* Gw=(\bar\pa\bar\pa^*+\bar\pa^*\bar\pa)Gw=w,
\]
because $w$ is $\bar\pa$-closed and because $\bar\pa\bar\pa^*$ commutes with $G$. Moreover, we have
\begin{align*}
&\int_{M}\|u\|^{2}_{h_{L}\otimes h_{E}}e^{-\Psi}dV_{g}=
\int_{M}\langle u, \bar\pa^* Gw\rangle_{h_{L}\otimes
h_{E}} e^{-\Psi} dV_{g}\\&=
\int_{M}\langle w, Gw\rangle_{g^{\ast}\otimes h_{L}\otimes
h_{E}} e^{-\Psi} dV_{g}\leq \frac 1C
\int_{M}\|w\|^{2}_{g^{\ast}\otimes h_{L}\otimes
h_{E}}e^{-\Psi} dV_{g}.
\end{align*}
Thus section $u$ defined above satisfies ~\eqref{2.4}.  For the rest of the paper, when we use Proposition~\ref{prop1}, we shall  fix this unique solution.
\end{rem}

Let $\eps=1/16$, and let $p'=[\eps(\log m)^2]$. Let $K_1$ be the upper bound of the curvature operator of $g$,
and let $K_2$ be the upper bound of the curvature operator of $h_E$. 
Throughout  the rest of the paper, we use $C$ to denote a positive constant, which may be different line by line,  depending only on $n$, $r$, $p',k_1,k_2$, and the injectivity radius $\delta$ of $M$.

Let $x_0\in M$ be a fixed point.
Let  $(z_1,\cdots,z_n)$  be holomorphic coordinates on 
an open neighborhood $U$ of $x_0$ in $M$. Define the function
$|z|$ by $$|z|=\sqrt{|z_1|^2+\cdots+|z_n|^2}$$  for $z\in U$. 
We further assume that $(z_1,\cdots,z_n)$ are the $K$-coordinates centered at $x_0$ and $e_L,\{e_\alpha\}$ are  the $K$-frames of order $p'$. We assume that $g_{i\bar j}, a$, and $h_{\alpha\bar\beta}$ are the local representatives of the metrics of $M,L$, and $E$, respectively.

Let $\Z_+^n$ be the set of $n$-tuple of integers
$(p_1,\cdots,p_n)$ such that $p_i\geq 0$ for  $(i=1,\cdots, n)$. Let
$P=(p_1,\cdots,p_n)$. Define
\begin{equation}\label{zp}
z^P=z_1^{p_1}\cdots z_n^{p_n}
\end{equation}
and
\[
p=p_1+\cdots+p_n.
\]

Let $\eta$ be  a smooth cut-off function
\begin{equation*}
\eta(t)=\begin{cases}1&\text{for $t<\frac{1}{2}$}\\0&\text{for $t\geq 1$}\end{cases}
\end{equation*} 
satisfying $0\leq -\eta'(t)\leq 4$ and $|\eta''(t)|\leq 8$.

The following lemma is a revised version of Tian's Lemma~\cite{t5}*{Lemma
1.2}.
\begin{lemma}\label{T}
For any $1\leq\alpha\leq r$ and    $P=(p_{1},\cdots,p_{n})\in
\mathbb{Z}_{+}^{n}$  such that $p'\geq p=p_{1}+\cdots +p_{n}$,
 there is a holomorphic global section $S_{P,\alpha,m}^{p'}$ of
$H^{0}(M,L^m\otimes E)$ satisfying
\begin{equation}\label{3+9}
\left|\int_{M}\|S_{P,\alpha,m}^{p'}\|_{h_L^m\otimes
h_E}^2 dV_g-1\right|\leq Ce^{-\frac{1}{4}(\log m)^{2}}.
\end{equation}
 $S^{p'}_{P,\alpha,m}$ can be decomposed as
\begin{equation}\label{0.51}
S_{P,\alpha,m}^{p'}=\widetilde{S}_{P,\alpha,m}-u_{P,\alpha,m}
\end{equation}
for $\widetilde{S}_{P,\alpha,m}$, $u_{P,\alpha,m}\in \mathcal C^\infty(M, L^m\otimes E)$. The support of $\widetilde{S}_{P,\alpha,m}$ is within $U$ and we can write
\begin{equation}
\widetilde{S}_{P,\alpha,m}(z)=\lambda_{P,\alpha}\eta\left(\frac{m|z|^2}{(\log
m)^2}\right)z^{P}e_L^m\otimes e_\alpha,\label{0.5}
\end{equation}
where
\begin{equation}
\lambda_{P,\alpha}^{-2}=\int_{|z|\leq\frac{\log
m}{\sqrt{m}}}|z^P|^2a^mh_{\alpha\bar \alpha}dV_{g}.\label{0.8}
\end{equation}
The smooth section $u_{P,\alpha,m}$, which  solves the equation
\begin{equation}\label{2345}
\bar\pa u_{P,\alpha,m}=\bar\pa \tilde S_{P,\alpha,m},
\end{equation}
satisfies 
\begin{equation}
u_{P,\alpha,m}(z)=O(|z|^{p'})\qquad\text{if $z\in U$},\label{0.6}
\end{equation}
and
\begin{equation}\int_{M}\|u_{P,\alpha,m}\|^2_{h_L^m\otimes h_E}dV_{g}\leq C\,{e^{-\frac{1}{4}(\log m)^{2}}}.\label{0.7}
\end{equation}
 \end{lemma}

\vspace{.1in}

Let $C(\mathbb R)$ be the space of real functions on $\mathbb R$ and let  the norm  of a function $f$ in $C(\mathbb R)$ be
\[
\|f\|=\sup_{x\in\mathbb R} |f(x) e^{-x}|,
\]
if the right hand side of the above is finite. 

For any positive integer $p$, the norm of the function  $x^p$ is $p^pe^{-p}$. The function  $x^pe^{-x}$  peaks at $x=p$. As $x\to\infty$, the function decays very fast. With this picture in mind, we justify  the following definition.

\begin{definition}\label{def31}  The sections $\{S^{p'}_{P,\alpha,m}\}$, which  look like $z^P$ near $x_0$,  are called peak sections. \end{definition}

\begin{proof}[Proof of Lemma~\ref{T}]
 Define the weight function
\[\Psi(z)=(2n+2p')\eta\left(\frac{m|z|^{2}}{(\log m)^{2}}\right)\log\left(\frac{m|z|^{2}}{(\log m)^{2}}\right).\]
A straightforward computation gives
\begin{equation}
\sqrt{-1}\partial\bar{\partial}\Psi\geq -\frac{8m(2n+2p')}{(\log
m)^{2}}\omega_{g}\geq -\frac{9}{10}m\,\omega_g\label{0.9}
\end{equation}
for $m$ large enough.
Let
\begin{equation}
\Theta=\Theta_{E}+I_{E}\otimes\partial\bar{\partial}
\Psi+\frac{2\pi}{\sqrt{-1}}I_{E}\otimes (Ric(h)+Ric(g)).
\end{equation}
Using
\eqref{0.9}, we verify that
\begin{align}\label{0.91}
\begin{split}
&\langle\Theta(\varphi),\varphi\rangle\geq
\left(\frac{1}{10}m-K_{1}-K_2\right)\|\varphi\|^2_{g^{\ast}\otimes h_L^m\otimes
h_{E}}
\\&\geq\frac{1}{20}m\|\varphi\|^{2}_{g^{\ast}\otimes h_L^m\otimes
h_E}
\end{split}
\end{align}
 for $\varphi\in \bigwedge^{0,1}(M,L^{m}\otimes
E)$.
 For $P\in \mathbb{Z}_{+}^{n}$, let
\begin{equation}\label{0.92}
w_{P,\alpha,m}=\bar{\partial}\tilde S_{P,\alpha,m}.
\end{equation}
 Since $\bar{\partial}
w_{P,\alpha,m}\equiv 0$  and
\[
\int_M\|w_{P,\alpha,m}\|_{g^*\otimes h_L^m\otimes h_E}^{2}e^{-\Psi}dV_{g}<+\infty,
\]
by Proposition \ref{prop1} and Remark~\ref{rem21},  there exists a smooth section $u_{P,\alpha,m}$ that solves ~\eqref{2345} such that 
\begin{equation}
\int_{M}\|u_{P,\alpha,m}\|^{2}_{h_L^m\otimes h_E}e^{-\Psi}dV_g
\leq\frac{20}{m}\int_{M}\|w_{P,\alpha,m}\|_{g^{\ast}\otimes
h_{L}^{m}\otimes h_{E}}^{2}e^{-\Psi}dV_{g}<\infty.\label{0.10}
\end{equation}
Since 
\[
\int_M\|u_{P,\alpha,m}\|^2 e^{-\Psi} dV_g<+\infty,
\]
we must  have $u_{P,\alpha,m}(x)=O(|z|^{p'})$.  This proves~\eqref{0.6}. In
order to verify~\eqref{0.7}, we need to estimate the right hand side
of~\eqref{0.10}. By the definition of
 function $\eta$, we have
\begin{align*}
\begin{split}
&\frac{20}m\int_M \|w_{P,\alpha,m}\|^2_{g^{\ast}\otimes h_L^m\otimes
h_E}e^{-\Psi} dV_g \\&\leq \frac{320}{m}\lambda^2_{P,\alpha}\int_{\frac{\log m}{\sqrt
{2m}}\leq\rr}\frac{m^{2}}{(\log
m)^{4}}z_{k}\bar{z}_{\ell}g^{\ell\bar{k}}|z^{P}|^{2} a^m h_{\alpha\bar
\alpha}e^{-\Psi}dV_g\\&\leq C\lambda^2_{P,\alpha}\frac{(\log
m)^{2(|P|-1)}}{m^{|P|}}\int_{\frac{\log m}{\sqrt {2m}}\leq\rr}a^m
dV_0,
\end{split}
\end{align*}
where $dV_0$ is the volume form of $\mathbb C^n$
\[
dV_0=\left(\frac{\sqrt{-1}}{2\pi}\right)^n dz_1\wedge d\bar z_1\wedge\cdots\wedge dz_n\wedge d\bar z_n.
\]
Since $e_L$ is
a $K$-frame, we have
\[
a=
e^{\log a}=e^{-|z|^2+O(|z|^4)}
\]
on $U$. Thus we have
\begin{align}\label{0.94-1}
\begin{split}
&\int_M\|u_{P,\alpha,m}\|_{h_{L}^{m}\otimes
h_{E}}^{2}e^{-\Psi}dV_g
\leq C\lambda^2_{P,\alpha} \frac{(\log m)^{2(|P|-1)}}{m^{|P|}} 
\int_{|z|\geq\frac{\log m}{\sqrt{2m}}} e^{-m|z|^2} dV_0
 \\&= C\lambda^2_{P,\alpha} \frac{(\log m)^{2(|P|-1+n)}}{m^{|P|+n}} e^{-\frac
12(\log m)^2}.
\end{split}
\end{align} 
On the other hand, by Lemma~\ref{lem321} below, for $m$ large enough
\[
\lambda_{P,\alpha}^{-2}\geq C\frac{P!}{m^{n+|P|}}
\]
and therefore
\begin{equation}\label{0.94}
\lambda_{P,\alpha}\leq Cm^{\frac{n+|P|}2}/\sqrt{P!}.
\end{equation}
Putting ~\eqref{0.94} and ~\eqref{0.94-1} together,  we obtain ~\eqref{0.7} for $m$ large enough.
Finally,  we have
\[
\left|\int_{M}\|\tilde S_{P,\alpha,m}\|_{h_{L}^{m}\otimes
h_{E}}^{2}dV_{g}-1\right|\leq
C\lambda_{P,\alpha}^2 \frac{(\log
m)^{2(|P|-1+n)}}{m^{|P|+n}}\int_{\frac{\log m}{\sqrt {2m}}\leq\rr}a^m
dV_g.
\]
Using the same estimate as above, we obtain
\[
\left|\int_{M}\|\tilde S_{P,\alpha,m}\|_{h_{L}^{m}\otimes
h_{E}}^{2}dV_{g}-1\right|\leq C e^{-\frac
13(\log m)^2}.
\]
for $m$ large enough. ~\eqref{3+9} then follows from ~\eqref{0.51} and all the above estimates.

\end{proof}

\begin{lemma}\label{lem321}
For any multi-index $P$ such that $|P|\leq\frac 14(\log m)^2$, we have
\begin{equation}\label{24-1}
\frac {P!}{m^{n+|P|}}(1-2^ne^{-\frac 14(\log m)^2})\leq \int_{|z|\leq\frac{\log m}{\sqrt m}}|z^P|^2 e^{-m|z|^2} dV_0\leq \frac{P!}{m^{n+|P|}}.
\end{equation}
\end{lemma}

\begin{proof}
We quote the following elementary identity:
\begin{equation}\label{ele}
\int_{\mathbb C^n}|z^P|^2|z|^{2q}e^{-m|z|^2} dV_0=\frac{(|P|+n-1+q)!P!}{(|P|+n-1)!m^{n+|P|+q}},
\end{equation}
where $q\geq 0$ and $P$ is a multi-index.

Using the above identity, the second inequality of ~\eqref{24-1}  follows. To prove the first inequality of~\eqref{24-1}, we need to show that 
\[
\int_{|z|\geq\frac{{\log  m}}{\sqrt m}}|z^P|^2 e^{-m|z|^2} dV_0\leq 2^n\frac{P!}{m^{n+|P|}}e^{-\frac 14(\log m)^2}.
\]

Rescaling  $z$ to $\sqrt mz$, the above inequality ~\eqref{24-1} becomes 
\begin{equation}\label{24-3}
\int_{|z|\geq{\log  m}}|z^P|^2 e^{-|z|^2} dV_0\leq 2^n P!e^{-\frac 14(\log m)^2}.
\end{equation}

We have
\[
 \int_{|z|\geq{\log  m}}|z^P|^2 e^{-|z|^2} dV_0\leq
 e^{-\frac 12(\log m)^2} \int_{|z|\geq{\log  m}}|z^P|^2 e^{-\frac 12|z|^2} dV_0,
 \]
 and by rescaling $z$ to $\sqrt 2z$, we obtain
 \begin{align*}
 & \int_{|z|\geq{\log  m}}|z^P|^2 e^{-|z|^2} dV_0\\
&\leq  e^{-\frac 12(\log m)^2} 2^{|P|+n} \int_{\mathbb C^n}|z^P|^2 e^{-|z|^2} dV_0=\frac{2^{|P|+n} P!}{m^{n+|P|}}e^{-\frac 12(\log m)^2}.
\end{align*}
The lemma then follows from the assumption
$|P|\leq\frac 14(\log m)^2$.

\end{proof}

\begin{definition}\label{def32}
For fixed $m,p'$, let $k$ be the smallest integer such that $|R(k)|\geq p'$. Let
\[
S_j=\lambda^{-1}_{R(j),\zeta(j)} \,S^{p'}_{R(j),\zeta(j),m}
\]
for $0\leq j\leq k$. It is apparent that $S_1,\cdots, S_k$ are linearly independent.
We extend $\{S_j\}_{j=1,\cdots,k}$ to $\{S_j\}_{j=1,\cdots,d}$ so that  $S_j=o(|z|^{p'})$ for $j>k$.
\end{definition}

\begin{cor}\label{cor31}  For any $x_0\in M$, the above basis $\{S_j\}$ is a strongly regular basis. Such a basis satisfies assumption (1) of Theorem~\ref{thm21}.
\end{cor}

\qed

\section{On  the generalization of a lemma of Ruan}\label{4}
W. D. Ruan discovered the almost  orthogonality between the peak sections and the sections of high orders.
For trivial vector bundle $E$, we quote Ruan's Lemma \cite{ru}*{Lemma 3.2} as follows.
\begin{lemma}[Ruan] \label{Ru}
Let $S_P=S_{P,m}^{p'}$ be a section constructed in Lemma \ref{T}. Define
\[\|S\|=\sqrt{\int_M|S|^2_{h_L^m\otimes g^\ast}dV_g}\] for $S\in H^0(M,L^m)$. Let $T$ be another section of $L^m$. Near $x_0$, $T=fe_L^m$ for a holomorphic function $f$. When we say $T$'s Taylor expansion at $x_0$, we mean the Taylor expansion of $f$ at $x_0$ under the coordinate system $(z_1,\cdots,z_n)$.
\begin{enumerate}
\item If $z^P$ is not in $T$'s Taylor expansion at $0$, then
\[(S_P,T)=O\left(\frac{1}{m}\right)\|T\|.\]
\item If $T$ contains terms $z^Q$ for $|Q|\geq |P|+\sigma$ in the Taylor expansion, then
\[
 (S_P,T)=O\left(\frac{1}{m^{1+\frac{\sigma}{2}}}\right)\|T\|.
 \]
 \end{enumerate}
\end{lemma}

In order to use Ruan's result,  we need to give the precise bound in the above  lemma.

We say that a smooth function $f$ on $U\subset \C^n$ is regular, if in its Taylor expansion, there are no $z^P, \bar z^P z^P\bar z_j, \bar z^P z_j$ terms. We begin with the following
 
\begin{lemma} 
Let $f$ be a holomorphic function on $U\subset \mathbb{C}^n$, let $a>0$ be a positive smooth function and $\xi$ a complex smooth function on $U$ such that $\xi(0)=1$. Assume that both $a$ and $\xi$ are regular functions. Define
\[\|f\|=\sqrt{\int_{|z|\leq\frac{\log m}{\sqrt{m}}}|f|^2 a^m dV_0}.\]
Assume that $|P|\leq\eps(\log m)^2$ for some $\eps>0$.
\begin{enumerate}
\item If $z^P$ is not in $f$'s Taylor expansion at $0$, then
\[\bigg|\int_{|z|\leq\frac{\log m}{\sqrt{m}}}z^P\bar{f} a^m\xi dV_0\bigg|\leq \frac{C\sqrt{P!}}{m^{1+\frac{n+|P|}{2}}}\|f\|.\]
\item If $f$ only contains terms $z^Q$ for $|Q|\geq |P|+\sigma$ and $\sigma>0$ in the Taylor expansion, then
\begin{equation}
 \bigg|\int_{|z|\leq\frac{\log m}{\sqrt{m}}}z^P\bar{f} a^m\xi dV_0\bigg|\leq \left(\frac{(\log m)^{10}}{m}\right)^{1+\frac{n+|P|+\sigma}{2}}\|f\|\label{r2}
 \end{equation}
 for $m>m(\sigma)$, where $m(\sigma)$ is a constant depending on $\sigma$ and the functions $a, \xi$.
 \item Moreover, if the metrics are analytic, then we have
 \begin{equation}
 \bigg|\int_{|z|\leq\frac{\log m}{\sqrt{m}}}z^P\bar{f} a^m\xi dV_0\bigg|\leq C^\sigma\left(\frac{(\log m)^{10}}{m}\right)^{1+\frac{n+|P|+\sigma}{2}}\|f\|\label{r3}
 \end{equation}
 for $m>m(\sigma)$.
 \end{enumerate}
\end{lemma}

\begin{proof}
  We shall omit the proof of (1), since it is similar to that of (2).

Let $\log a+|z|^2=\zeta$. 
Since $a$ is a regular function, $\zeta=O(|z|^4)$. It follows that 
\begin{equation}\label{qsq}
 |m\zeta|\leq \frac{1}{{m^{3/4}}}
 \end{equation}
for   $|z|\leq\frac{\log m}{\sqrt{m}}$    and $m$ is sufficiently large.  As a result, the functions $a^m$ and $e^{-m|z|^2}$ are mutually bounded by a constant independent to $m$. In particular, \[
\|f\|\leq C \sqrt{\int_{|z|\leq\frac{\log m}{\sqrt{m}}}|f|^2 e^{-m|z|^2} dV_0}
\]
for some constant $C>0$. We shall use this fact below repeatedly without further notice.

Let  $\zeta=\zeta_1+\zeta_2$ such that $\zeta_1$ is the Taylor's polynomial of $\zeta$ of order  $2\sigma+1$, and 
  let $\xi=\xi_1+\xi_2$ such that $\xi_1$ is the Taylor's polynomial of $\xi$ of order  $2\sigma+1$. 
  Then for $m$ large enough, we have
 \begin{equation}\label{qsq-0}
 |e^{m\zeta}\xi-\sum_{k=0}^{2\sigma+1}\frac{m^k}{k!}\zeta_1^k\xi_1|
 \leq \frac{1}{m^{1+\frac\sigma 2}}
 \end{equation}
 for  $|z|\leq\frac{\log m}{\sqrt{m}}$.
Therefore using the Cauchy inequality, we get
  \begin{align}\label{qsq-1}
  \begin{split}
  &\bigg|\int_{|z|\leq\frac{\log m}{\sqrt{m}}}z^P\bar{f} e^{-m|z|^2}     
  (e^{m\zeta}\xi-\sum_{k=0}^{2\sigma+1}\frac{m^k}{k!}\zeta_1^k\xi_1)   dV_0\bigg|\\
  &\leq \frac{C\sqrt{P!}}{m^{1+\frac{n+|P|+\sigma}{2}}}\|f\|\leq  
   \left(\frac{(\log m)^2}{m}\right)^{1+\frac{n+|P|+\sigma}{2}}\cdot\|f\|
  \end{split}
  \end{align}
  by~\eqref{ele} and the fact that $|P|\leq \eps (\log m)^2$.

   Note that $\zeta_1^k\xi_1$ is a  polynomial in $z$ and $\bar{z}$. Let
  \begin{equation}
  \zeta_1^k\xi_1=\sum_{IJ}\zeta_{IJ}z^I\bar{z}^J.\label{zeta}
  \end{equation}
The coefficients
\begin{equation}\label{qsq-2}
|\zeta_{IJ}|=\frac{1}{I!J!}\left|\frac{\pa^{|I|+|J|}(\zeta_1^k\xi_1)}{\pa z^I\pa\bar z^J}(0)\right|=\frac{1}{I!J!}\left|\frac{\pa^{|I|+|J|}(\zeta_1^k\xi)}{\pa z^I\pa\bar z^J}(0)\right|\leq C (\sigma)
\end{equation}
for $C(\sigma)>0$ and for $m$ large enough. 
  
   If $|I|-|J|<\sigma$, then by the assumption on $f$ and by symmetry,
  \begin{equation}\label{symm}
  \int_{|z|\leq\frac{\log m}{\sqrt{m}}}z^P\bar{f} e^{-m|z|^2} m^{k} \zeta_{IJ}z^I\bar{z}^JdV_0=0.
  \end{equation}
  On the other hand, under the $K$-coordinates and $K$-frames, in the expansion of $\zeta$, there is no $z^P$ or $z^P\bar{z_j}$ terms. Thus in \eqref{zeta}, we must have $|J|\geq 2k$. If $|I|-|J|\geq \sigma$ and $|J|\geq 2$, then
  \[|I|+|J|-2k\geq 2+\sigma.\]
  Hence that
  \begin{equation}\label{qsq-3}
  \left|\sum_{|I|+|J|-2k\geq 2+\sigma}\zeta_{IJ}z^I z^{\bar J}\right|\leq C(\sigma)|z|^{2k+2+\sigma}.
  \end{equation}
  
Using ~\eqref{qsq-2},~\eqref{symm}, ~\eqref{ele} and the above estimate, we have
  \begin{align*}
  &
  \bigg|\int_{|z|\leq\frac{\log m}{\sqrt{m}}}z^P\bar{f} e^{-m|z|^2} m^{k}\zeta_1^k\xi_1 dV_0\bigg|\\
  &\leq
 m^k\|f\|\cdot \sqrt{\int_{|z|\leq\frac{\log m}{\sqrt{m}}}|z^P|^2\left|\sum_{|I|+|J|-2k\geq 2+\sigma}\zeta_{IJ}z^I z^{\bar J}\right|^2e^{-m|z|^2} dV_0}\\
 & \leq C(\sigma) (\log m)^{-1}\left(\frac{(\log m)^{10}}{m}\right)^{1+\frac{n+|P|+\sigma}{2}}\cdot\|f\|.
  \end{align*}
  Combining the above inequality and ~\eqref{qsq-0}, ~\eqref{qsq-1}, we proved that 
  \begin{equation}
 \bigg|\int_{|z|\leq\frac{\log m}{\sqrt{m}}}z^P\bar{f} a^m\xi dV_0\bigg|\leq C(\sigma)(\log m)^{-1}\left(\frac{(\log m)^{10}}{m}\right)^{1+\frac{n+|P|+\sigma}{2}}\cdot \|f\|.
 \end{equation}
 Since $m$ is sufficiently large, the above inequality implies ~\eqref{r2}.
 
 The proof of (3) is similar. We notice that in the analytic case, ~\eqref{qsq-0} becomes
 \[
  |e^{m\zeta}\xi-\sum_{k=0}^{2\sigma+1}\frac{m^k}{k!}\zeta_1^k\xi_1|
 \leq \frac{C^\sigma}{m^{1+\frac\sigma 2}}
 \]
 for some constant $C>0$
 and ~\eqref{qsq-1} becomes
 \begin{align*}
 &\bigg|\int_{|z|\leq\frac{\log m}{\sqrt{m}}}z^P\bar{f} e^{-m|z|^2}     
  (e^{m\zeta}\xi-\sum_{k=0}^{2\sigma+1}\frac{m^k}{k!}\zeta_1^k\xi_1)   dV_0\bigg|\\&
 \leq  
   C^{\sigma}\left(\frac{(\log m)^2}{m}\right)^{1+\frac{n+|P|+\sigma}{2}}\cdot\|f\|
\end{align*}
for a possibly larger constant $C$.

Therefore, we have
 \[
 \left|\sum_{|I|+|J|-2k\geq 2+\sigma}\zeta_{IJ}z^I z^{\bar J}\right|\leq C^\sigma|z|^{2k+2+\sigma}\]
 for some constant $C$. (3) follows from the above analytic version of ~\eqref{qsq-0},~\eqref{qsq-1}, and ~\eqref{qsq-3}.
This completes the proof of the lemma.

\end{proof}

By taking  $\xi=(\det g_{i\bar{j}}) h_{\alpha\bar{\beta}} $, we proved the following results, which extend the estimates up to $[\eps(\log m)^2]$.

\begin{theorem}\label{ruan-mod}
Let $S_{P,\alpha}=S_{P,\alpha,m}^{p'}$ for $|P|\leq p'$, $1\leq \alpha \leq r$ be the peak sections defined in Lemma~\ref{T} and let $T$ be a section with vanishing order at $x_0$ at least $p'$. Then for $m>m(\sigma)$, we have
\begin{enumerate}
\item
\[
|(S_{P,\alpha},T)|\leq\frac{C}{m^{3/2}}\|T\|.
\]
\item Moreover, 
\[
|(S_{P,\alpha},T)|\leq \frac{C}{m^{1+\frac 12({p'-|P|})}}\|T\|.
\]
\item If the metrics are analytic, then
\[
|(S_{P,\alpha},T)|\leq \frac{C^{\sigma}}{m^{1+\frac 12({p'-|P|})}}\|T\|.
\]
\end{enumerate}
\end{theorem}

\qed

Taking a special case, we can prove that

\begin{cor} \label{cor41} The basis $\{S_j\}$ is an almost orthonormal basis and assumption (2) of Theorem~\ref{thm21} is satisfied.
\end{cor}

\qed

Using the Taylor's expansion, we have

\begin{theorem}\label{thm4.2}
For fixed $i,j$ such that $|R(i)|+|R(j)|\leq\mu$ and $\alpha=\zeta(i),\beta=\zeta(j)$, there is an asymptotic expansion 
\[
(S_i,S_j)\sim \frac{1}{m^{n+\frac 12(|R(i)|+|R(j)|)}}\left(A_{0ij}+\frac{A_{1ij}}{m}+\cdots\right)
\]
for $r\times r$ matrices $A_i$ in the sense that for any $\mu>0$, 
\begin{enumerate}
\item  If the metrics are smooth, then for any $s\leq\eps(\log m)^2$, we have
\[
\left|(S_i,S_j)-\frac{1}{m^{n+\frac 12(|R(i)|+|R(j)|)}}\left(A_{0ij}+\frac{A_{1ij}}{m}+\cdots+\frac{A_{sij}}{m^s}\right)\right|\leq \frac{C}{m^{s+1}}.
\]
\item If the metrics are real analytic, then for $s=[\eps(\log m)^2]$, we have
\[
\left|(S_i,S_j)-\frac{1}{m^{\frac{n+|R(i)|+|R(j)|}{2}}}
\left(A_{0ij}+\frac{A_{1{ij}}}{m}+\cdots
\frac{A_{s{ij}}}{m^s}\right)\right|\leq \frac{C^{s+1}}{m^{s+1}},
\]
where
\[
\|A_{\xi ij}\|_{\mathcal C^\mu}\leq C^\xi
\]
for $\xi\geq 0$.
\end{enumerate}
\end{theorem}

\begin{proof}
In~\cite{lu10}*{Proposition 2.1}, (1) was proved for  the trivial line bundle $E$. The proof of the general case is the same so we omit the proof.  
If the metrics are real analytic, the Taylor expansions of all the metric matrices are convergent. Therefore using the same method as in (1), we prove the second part of the theorem.

\end{proof}

\begin{cor}\label{cor51}
The assumptions (3), (4) of Theorem~\ref{thm21} are valid.
\end{cor}

\begin{proof}[Proofs  of Theorem~\ref{main-2} and Theorem~\ref{thm1.2}.] These two theorems follow from Theorem~\ref{thm21}, Corollary~\ref{cor31}, Corollary~\ref{cor41}, and Corollary~\ref{cor51}.

\end{proof}

\begin{proof}[Proof of Corollary~\ref{cor11}]
This is the same as the proof of ~\eqref{cor11-1}.

\end{proof}

\appendix

\section{$K$-coordinates and $K$-frames}\label{sec:6}
Let $\mu$ be a nonnegative integer. 
It is well-known that other than $\mu=0$,  the $\mathcal C^\mu$-norm  on the space of smooth functions depends on the choice of local coordinate systems. For two different atlases of the manifold, the two different $\mathcal C^\mu$-norms are equivalent (mutually bounded).
Therefore the underlying topology is intrinsically defined by these norms. 

In the proof of the main results of this paper, we need to treat {\it uncountably many} local coordinate systems. Therefore, it is necessary to look into the details of the definition of $\mathcal C^\mu$-norms.

Let $P$ be a multiple index: $P=(p_1,\cdots,p_n)$, where $p_1,\cdots,p_n$ are nonnegative numbers. Let $|P|=p_1+\cdots+p_n$, and let $P!=p_1!\cdots p_n!$. Define
\[
z^P=z_1^{p_1}\cdots z_n^{p_n}.
\]
Let $f$ be a smooth  function on an open set $U$ of $M$ with holomorphic local coordinates $z=(z_1,\cdots,z_n)$. Define
\[
|D^k f|=\sum_{|P|+|Q|=k}\frac{k!}{P!Q!}| D^{P,Q} f|
\]
for $k\geq 0$,
where
\[
P=(p_1,\cdots,p_n),\qquad Q=(q_1,\cdots,q_n),
\]
and $D$ is  the differential operator
\[
D^{P,Q}f=\frac{\pa^{|P|+|Q|}f}{\pa z_1^{p_1}\cdots \pa z_n^{p_n}\pa \bar z_{1}^{q_1}\cdots \pa\bar z_n^{q_n}}.
\]
$|D^kf|$ defines a nonnegative function on $U$.
The $\mathcal C^\mu$-norm of $f$ on $U$ is defined by
\[
\|f\|_{\mathcal C^\mu}=\max_{k\leq \mu}\sup_{x\in U} |D^kf|(x).
\]
The definition of $|D^k f|$ depends on the local coordinates. We use the notation $|D^k f|_z$ to denote such a dependence. Let $w=(w_1,\cdots,w_n)$ be another local coordinates on $U$, then there is a constant $C>1$, depending on $z,w$, such that
\[
C^{-1}|D^kf|_z\leq |D^kf|_w\leq C|D^kf|_z.
\]
\\

Let $M$ be an $n$-dimensional algebraic manifold with a positive
Hermitian line bundle $(L,h_L)\rightarrow M$.
Let $(E,h_E)$ be a Hermitian vector bundle of rank $r$ over $M$.
Suppose that the
K\"ahler form $\omega_g$ of the \ka metric $g$ is defined by the
curvature ${\rm Ric}(h_L)$ of $h_L$. That is, under local coordinates $(z_1,\cdots,z_n)$ at a fixed point $x_0$, we have
\[
\omega_g=-\bb\sum_{\alpha,\beta=1}^n \frac{\pa^2}{ \pa
z_\alpha\pa\bar z_\beta}\log a\, dz_\alpha\wedge d\bar
z_\beta=\bb\sum_{\alpha,\beta=1}^n g_{\alpha\bar\beta}
dz_\alpha\wedge d\bar z_\beta,
\]
where $a$ is the local representation of  the Hermitian metric
$h_L$.

\begin{definition}
 Let $p>0$ be any positive integer. Let $x_0\in M$ be a point. Let 
 $(z_1,\cdots,z_n)$ be a holomorphic coordinate system centered at $x_0$. Let 
$(g_{\alpha\bar\beta})$ be the  \ka metric matrix. If it satisfies
\begin{align}\label{2-04}
\begin{split}
& g_{\alpha\bar\beta}(x_0)=\delta_{\alpha\beta};\\
&\frac{\pa^{p_1+\cdots+p_n}g_{\alpha\bar\beta}}{\pa
z_1^{p_1}\cdots\pa z_n^{p_n}} (x_0)=0
\end{split}
\end{align}
for $\alpha, \beta=1,\cdots,n$ and any nonnegative integers $(p_1,\cdots,p_n)$ with
$p>p_1+\cdots+p_n\neq 0$. Then we call the coordinate system a $K$-coordinate system of order $p$.
\end{definition}

The existence  of 
$K$-coordinate system was known  to the string theorists in the 1980s. However, the result had been known to Bochner long time ago.
We refer to~\cite{bo} for the proof of the existence of $K$-coordinates.

If the metric is analytic, then we can take $p$ to be $+\infty$.  In this case, the $K$-coordinate system is unique up to an affine transformation.

Similar to the above, we make the following definition

\begin{definition}
Let  $e_L$  be a local holomorphic frame of $L$ at $x_0$. If 
for $p>0$, the local representation function $a$ of
the Hermitian metric $h_L$ satisfies
\begin{equation}\label{2-05}
a(x_0)=1, \frac{\pa^{p_1+\cdots+p_n}a}{\pa z_1^{p_1}\cdots\pa
z_n^{p_n}}(x_0)=0
\end{equation}
for any nonnegative integers $(p_1,\cdots,p_n)$ with
$p>p_1+\cdots+p_n\neq 0$. Then we call $e_L$ is a $K$-frame of order $p$. If $a$ is analytic, then again we can take $p=+\infty$.
\end{definition}

The following lemma is similar to the results above. We sketch the proof here.

\begin{lemma} For any $p>0$,
it is possible to choose a local holomorphic frame $\{e_1,\cdots,e_r\}$ of order $p$ on $U$ such that
\begin{align}\label{2.6}
\begin{split}
& h_{i\bar j}(x_0)=\delta_{ij};\\
&\frac{\pa^{p_1+\cdots+p_n}h_{i\bar j}}{\pa
z_1^{p_1}\cdots\pa z_n^{p_n}} (x_0)=0
\end{split}
\end{align}
for any nonnegative integers $(p_1,\cdots,p_n)$ with $p>p_1+\cdots+p_n\neq 0$,
where $h_{i\bar j}=\langle e_i,e_j\rangle$. If the metric is analytic, then we can take $p=+\infty$.
Moreover, any derivative of $h_{i\bar j}$ at $x_0$ can be represented as a polynomial of curvatures  of both  $E$ and $M$ and their derivatives.
\end{lemma}

\begin{proof} We first choose a holomorphic frame $\{e_1,\cdots,e_r\}$ such that
$h_{i \bar j}(x_0)=\delta_{i j}$. Let the matrix $H=(h_{i\bar j})$. Let the Taylor expansion of $H$ be
\[
H\sim I+A+\bar A+B,
\]
where $I$ is the identity matrix; $A$ is the holomorphic part of the Taylor expansion; $\bar A$ is the complex conjugate of $A$;
and  $B$ is the mixed part, that is, the entries of $B$ are composed of both $z$'s and $\bar z$'s.
~\footnote{If $H$ is analytic, the Taylor expansion must be convergent.
If $H$ is smooth, the expansion is understood as {\it formal}:
it doesn't have to  converge, and even if it does, it doesn't have to  converge to the matrix-valued function $H$.}

Since $\bar H^T=H$, we must have $A=A^T$.
Let $p>2$ be an integer. Let $A_p$ be the first $p$ terms in the formal series $A$.
 Define a new frame $\{f_1,\cdots,f_r\}$ such that
\[
e_i=f_i+(A_p)_{ij}f_j.
\]
It is not hard to see that under the new frame, the metric matrix is $\tilde H=(I+A_p)^{-1}H(I+\bar A_p)^{-1}$.
A straightforward computation shows that in the Taylor expansion of $\tilde H$,
there are no holomorphic or anti-holomorphic parts up to order $p$.

Finally, we have the formula
\[
\frac{\pa^2 h_{i\bar j}}{\pa z_\alpha\pa\bar z_\beta}=(\Theta_E)_{i\bar j\alpha\bar \beta}
+h^{k\bar l}\frac{\pa h_{i\bar l}}{\pa z_\alpha}\cdot\frac{\pa h_{k\bar j}}{\pa \bar z_\beta}.
\]

By the above equation and by induction, all derivatives of $(h_{i\bar j})$ at $x_0$ can be expressed as polynomials of $\Theta_E$,
the curvature of $M$, and their covariant derivatives.

\end{proof}

Because of the above lemma, we will call the  frames $e_L$ and $\{e_j\}$ $K$-frames,
 that is, the  local holomorphic coordinates, the local frames of $L$,  and the local frame of $E$ will satisfy ~\eqref{2-04},~\eqref{2-05}, and~\eqref{2.6}.

As we have discussed before,  $K$-coordinates and $K$-frames are not unique.
However, in what follows, we will  write out explicitly a smooth family of $K$-frames and $K$-coordinates near any given point.

Let
\[
a(z)=\langle e_L,e_L\rangle,\qquad H(z)=(h_{i\bar j}),\qquad G(z)=(g_{i\bar j}).
\]
Let $P$ be the multiple index: $P=(p_1,\cdots,p_n)$. We define $|P|=\sum p_i$;
$P!=p_1!\cdots p_n!$; $z^P=z_1^{p_1}\cdots z_n^{p_n}$; and
\[
f^{(P)}(z)=\frac{\pa^P f}{\pa z^P}(z)
\]
for a smooth function $f$. Using these notations,
~\eqref{2-04},~\eqref{2-05}, and~\eqref{2.6} can be written as
\[
g_{\alpha\bar\beta}(x_0)=\delta_{\alpha\beta}, a(x_0)=1, h_{i\bar j}(x_0)=\delta_{ij},
\]
and
\[
(g_{\alpha\bar\beta})^{(P)}(x_0)=0, a^{(P)}(x_0)=0, (h_{i\bar j})^{(P)}(x_0)=0
\]
for $0<|P|<p$.

 For any $p>2$ and any $t\in U$  with $|t|$ very small, we define
\begin{align}
\begin{split}
& a_{ij}(t)=\left.\frac{\pa^2\log a(z)}{\pa z_i\pa \bar z_j}\right|_{z=t},\\
& b_1^p(z,t)=\sum_{1\leq |P|\neq 0}\frac{(\log a)^{(P)}(t)}{P!}(z-t)^P,\\
&b_2^p(z,t)=\sum_{j=1}^n\sum_{1<|P|\leq p}\frac{1}{P!}\frac{\bar \pa(\log a)^{(P)}(t)}{\pa \bar t_j}\,(z-t)^P(\bar{z_j-t_j}).
\end{split}
\end{align}
We write
\begin{align}\label{plk}
\begin{split}
&
\log a(z)=\log a(t)+a_{ij}(t)(z_i-t_i)(\bar{z_j-t_j})\\
&+b_1^p(z,t)+\bar{b_1^p(z,t)}+b_2^p(z,t)+\bar{b_2^p(z,t)}+c^p(z,t).
\end{split}
\end{align}
By the definition of $c^p$, we have
\[
\left.\frac{\pa^P c^p}{\pa z^P}\right|_{z=t}=0,\quad
\left.\frac{\pa}{\pa \bar z_j}\frac{\pa^P c^p}{\pa z^P}\right|_{z=t}=0
\]
for $|P|\leq p$ and $1\leq j\leq n$.
Let $G^p(z,t), R^p(z,t)$ be the matrix-valued functions such that
\[
G^p(z,t)_{ij}=-\frac{\pa^2 b_2^p(z,t)}{\pa z_i\pa \bar z_j},\quad R^p(z,t)_{ij}=-\frac{\pa^2 c^p(z,t)}{\pa z_i\pa \bar z_j}.
\]
Then we have
\begin{equation}\label{edr}
G(z)=G(t)+G^p(z,t)+\bar{G^p(z,t)}^T+R^p(z,t)
\end{equation}
with
\begin{equation}\label{edr-1}
\left.\frac{\pa^P R^p}{\pa z^P}\right|_{z=t}=0
\end{equation}
for $0\leq |P|\leq p-1$.

Let $D(t)$  be a smooth Hermitian matrix-valued function such that $D(0)=I$ and $D(t)^2=G(t)$.
 Such a function exists. For example, if
 $\sum a_j x^j$ is the Taylor expansion of the function $\sqrt{1+x}$, then we can define $D(t)=\sum a_j (G(t)-I)^j$.

We use the following matrix notations: $z'=(z_1',\cdots,z_n')$, $z=(z_1,\cdots,z_n)$, $t=(t_1,\cdots,t_n)$ and
\[
\frac{\pa b_2^p}{\pa \bar z}=(\frac{\pa b_2^p}{\pa\bar  z_1},\cdots,\frac{\pa b_2^p}{\pa\bar  z_n}).
\]
We define
\begin{equation}\label{0.93}
z'=zD(t)-tD(t)-\frac{\pa b_2^p(z,t)}{\pa \bar z} D(t)^{-1}.
\end{equation}
Then $z'$ is a holomorphic coordinate system centered at $t$.
Differentiating the above, we have
\[
dz'=dz(D(t)+G^p(z,t)D(t)^{-1}).
\]
Under the new coordinates, the metric-matrix is
\[
P(z,t) ^{-1} G(z)(P(z,t)^{-1})^*,
\]
where $P(z,t)=D(t)+G^p(z,t)D(t)^{-1}$.
By~\eqref{edr}, the above matrix is equal to
\[
I+P(z,t)^{-1}(R^p(z,t)-G^p(z,t)G^{-1}(t)\bar {G^p(z,t)}^T)(\bar{P(z,t)}^{-1})^T.
\]
If $z=t$, then the above matrix is the identity matrix. Moreover,
since $P(z,t)$ is holomorphic with respect to $z$ and $G^p(t,t)=0$, from~\eqref{edr-1}, we have
\[
\left.\frac{\pa^P}{\pa z^P}\right|_{z=t}P(z,t)^{-1}(R^p(z,t)-G^p(z,t)G^{-1}(t)\bar {G^p(z,t)})(\bar{P(z,t)}^{-1})^T=0
\]
for any $|P|<p-1$.  By the chain rule, we conclude that
$z'$ is a $K$-coordinate system for any $t$.

Using a similar way, we can construct $K$-frames as well.

We let
\begin{align}
\begin{split}
&\xi^p(z,t)=\sum_{|P|\neq 0}^{|P|\leq p}\frac{ a^{(P)}(t)}{P!}(z-t)^P;\\
&B^p(z,t)=\sum_{|P|\neq 0}^{|P|\leq p}\frac{ H^{(P)}(t)}{P!}(z-t)^P,
\end{split}
\end{align}
and we write
\begin{align}
\begin{split}
&a(z)=a(t)+\xi^p(z,t)+\bar{\xi^p(z,t)}+\eta^p(z,t);\\
&H(z)=H(t)+B^p(z,t)+\bar{B^p(z,t)^T}+C^p(z,t).
\end{split}
\end{align}

As before, we let $K(t)$ be a smooth Hermitian matrix-valued function such that $K(t)^2=H(t)$.
We can rewrite
\begin{align}
& a(z)=\left(\sqrt{a(t)}+\frac{\xi^p(z,t)}{\sqrt{a(t)}}\right)\left(\sqrt{a(t)}+\frac{\bar{\xi^p(z,t)}}{\sqrt{a(t)}}\right)+\tilde\eta^p(z,t);\label{0.92-1}\\
&H(z)=(K(t)+B^p(z,t)K(t)^{-1})(K(t)+B^p(z,t)K(t)^{-1})^*+\tilde C^p(z,t),\label{0.92-2}
\end{align}
where
\begin{align*}
&\tilde\eta^p(z,t)=\eta^p(z,t)-\frac{|\xi^p(z,t)|^2}{a(t)};\\
&\tilde C^p(z,t)=C^p(z,t)-B^p(z,t)H^{-1}(t)\bar{B^p(z,t)}^T.
\end{align*}

As above, we know that
\begin{align}
&\left.\frac{\pa^P}{\pa z^P}\right|_{z=t}\tilde\eta^p(z,t)=0;\label{2.13}\\
&\left.\frac{\pa^P}{\pa z^P}\right|_{z=t}\tilde C^p(z,t)=0\label{2.14}
\end{align}
for $|P|<p$.

Define
\begin{align}
& e_L(t)=\left(\sqrt{a(t)}+\frac{\xi^p(z,t)}{\sqrt{a(t)}}\right)^{-1} e_L;\label{0.93-1}\\
& e(t)=e(\left(K(t)+B^p(z,t)K(t)^{-1}\right)^{-1})^T,\label{0.93-2}
\end{align}
where $e=(e_1,\cdots,e_r)$ and $e(t)=(e_1(t),\cdots,e_r(t))$.

Then by ~\eqref{0.92-1}, we have
\[
\langle e_L(t), e_L(t)\rangle=1+\tilde\eta^p(z,t)\left|\sqrt{a(t)}+\frac{\xi^p(z,t)}{\sqrt{a(t)}}\right|^{-2}.
\]
By~\eqref{2.13}, we know that in the Taylor expansion of $\langle e_L(t), e_L(t)\rangle$ at $t$, there are no holomorphic part up to the order $p$.
This proves that $e_L(t)$ is a $K$-frame for any $t$. Similarly, let
\[
Q(z,t)=K(t)+B^p(z,t)K(t)^{-1}.
\]
Then by~\eqref{0.92-2}, we have
\[
e(t)^T\bar{e(t)}=I+Q(z,t)^{-1}\tilde C^p(z,t)(\bar{Q(z,t)}^{-1})^T.
\]
By ~\eqref{2.14}, $e(t)$ is a $K$-frame for any $t$.

In summary, we prove

\begin{lemma} [The first stability lemma]
For $|t|$ small,  \eqref{0.93}, \eqref{0.93-1}, and ~\eqref{0.93-2} define a smooth family of  $K$-frames and $K$-coordinates. Moreover, for any $t$, the $K$-frames and $K$-coordinates are defined on the set
\[
|z'|<\delta/3,
\]
where $\delta$ is the injectivity radius of $M$ and $(z'_1,\cdots,z_n')$ are the $K$-coordinate system at $t$.
\end{lemma}

\qed

\begin{cor}\label{corA1}
Let $f$ be a smooth function on $M$. Assume that at any $x_0\in M$ and for any $K$-coordinates $(z_1,\cdots,z_n)$ at $x_0$, 
\[
\sup_{x_0,k\leq\mu}|D^k f|(x_0)\leq 1.
\]
Then there is a constant $C(\mu)$, depending on $\mu$, such that
\[
\|f\|_{\mathcal C^\mu}\leq C(\mu).
\]
\end{cor}

\qed

By the above lemma and the continuity, the following must be true:
for any $+\infty> p>0$ and  $x\in M$, there exists a $\rho=\rho_{x,p}>0$ such that
\begin{enumerate}
\item for each $|t|<\rho$, the smooth $K$-frames and $K$ coordinates exist and are at least of size $\rho$. That is, $\{|z'|<\rho\}$ is  contained in the $K$-coordinate system;
\item Define
\begin{align*}
&B_{1,k}(t)=\sum_{\ell\leq k}\frac{1}{\ell!}\|D^\ell\log a(t)\|(\delta/3)^\ell;\\
&B_{2,k}(t)=\max_{i,j}\sum_{\ell\leq k}\frac{1}{\ell!}\|D^\ell h_{Ei\bar j}(t)\|(\delta/3)^\ell;\\
&
A_k(t)=\sum_{k_1+k_2=k}B_{1,k_1}(t)\cdot
B_{2,k_2}(t),
\end{align*}
where $a(t)$ and $h_{Ei\bar{j}}(t)$ are the metric representation under the $K$-frames and $K$-coordinates. Then $A_k(t)\leq 2 A_k(0)$ for any $p>k\geq 0$.
\item  Let $U_t=\{|z'|<\rho\}$ be the $K$-coordinate neighborhood. Then on $U_t$, we have
\begin{equation}\label{nmp}
\log a(t)>\frac 12,\qquad g_{ij}(t)>\frac 12\delta_{ij},\qquad h_{i\bar j}\geq\frac 12\delta_{ij}.
\end{equation}
\end{enumerate}
Since $M$ is compact, finitely many of the above neighborhoods cover $M$. Therefore, we are able to define the norms $A_k$ as in the last section.

We have the following
\begin{lemma}[The second stability lemma] For any $+\infty\geq p>0$ and
for any point $x\in M$, there exist $K$-frames and $K$-coordinates of order $p$ such that
\begin{enumerate}
\item $A_k(t)\leq 2A_k(0)$, where $t\in U_t$;
\item On $U_t$,
~\eqref{nmp} is valid.
\end{enumerate}
\end{lemma}

\begin{proof}
The lemma essentially follows from continuity and compactness. 

\end{proof}

\end{document}